\pdfoutput=1
\RequirePackage{ifpdf}
\ifpdf % We are running pdfTeX in pdf mode
\documentclass[pdftex]{sigma}
\else
\documentclass{sigma}
\fi

\numberwithin{equation}{section}

\newtheorem{Theorem}{Theorem}[section]
\newtheorem{Conjecture}[Theorem]{Conjecture}
\newtheorem{Problem}[Theorem]{Problem}
\newtheorem{Question}[Theorem]{Question}
 { \theoremstyle{definition}
\newtheorem{Definition}[Theorem]{Definition}
 }

\begin{document}
\allowdisplaybreaks

\newcommand{\arXivNumber}{2007.04460}

\renewcommand{\thefootnote}{}

\renewcommand{\PaperNumber}{095}

\FirstPageHeading

\ShortArticleName{Covariant vs Contravariant Methods in Differential Geometry}

\ArticleName{Covariant vs Contravariant Methods\\ in Differential Geometry\footnote{This paper is a~contribution to the Special Issue on Scalar and Ricci Curvature in honor of Misha Gromov on his 75th Birthday. The full collection is available at \href{https://www.emis.de/journals/SIGMA/Gromov.html}{https://www.emis.de/journals/SIGMA/Gromov.html}}}

\Author{Maung MIN-OO}

\AuthorNameForHeading{M.~Min-Oo}

\Address{McMaster University, Hamilton, Ontario, Canada}
\Email{\href{mailto:minoo@mcmaster.ca}{minoo@mcmaster.ca}}

\ArticleDates{Received July 14, 2020, in final form September 17, 2020; Published online September 30, 2020}

\Abstract{This is a short essay about some fundamental results on scalar curvature and the two key methods that are used to establish them.}


\Classification{53C20; 53C21; 53C24; 53C27}

\renewcommand{\thefootnote}{\arabic{footnote}}
\setcounter{footnote}{0}

\section{Introduction}

\looseness=1 The first time that I got intrigued by the title of this essay was about 50 years ago when I learned the Bonnet--Meyers theorem which gives a sharp upper bound for the diameter in terms of a~positive lower bound for the Ricci curvature. This implies finiteness of the fundamental group and hence the vanishing of the first Betti number. In particular, the torus cannot carry a metric of positive Ricci curvature. The vanishing of the first de Rham cohomology group for closed manifolds with a metric of positive Ricci curvature can also be proved by using harmonic one-forms, as was first done by Bochner. Bochner's vanishing theorem is, in a sense, weaker than the diameter estimate, but it uses a completely different approach, which I call the {\it contravariant} method, in contrast to the classical approach using variations of geodesics and Jacobi fields, which I call the {\it covariant} method. The second time was about ten years after that when I learned about the proofs of the positive mass conjecture in general relativity. Schoen and Yau~\cite{SY3} used the covariant method, based on minimal surfaces whereas Witten~\cite{W1} used the contravariant approach with spinors and the Dirac operator. Witten's approach seems a bit weaker since one needs to assume that the manifold is spin, but that is not inappropriate for physical reasons (supersymmetry). Besides, Witten's proof reminded me of an earlier seminal paper by Gromov and Lawson~\cite{GL1} on why the torus cannot carry a metric of positive scalar curvature. Their proof was also contravariant using twisted Dirac operators and the Atiyah--Singer index theorem. In dimension~$3$ this result was first proved by Schoen and Yau~\cite{SY1} using minimal surfaces (covariant methods).

I first heard of Gromov in 1971 when I was still a student. I had the privilege of meeting him for the first time a few years later at the Arbeitstagung in Bonn. He explained to me a number of different things. I listened carefully and tried to understand as much as I could (and that hasn't changed during the last 45 years, every time I had a chance to speak with him!). He has been a truly inspirational figure in my mathematical career.

\section{Jacobi fields vs Dirac spinors}

\subsection{Covariant methods}

By this I mean methods using mainly geodesics, minimal surfaces, first and second variation formulas, cut-locus estimates, comparison theorems, etc. It is the standard approach that I~learned during the 70s and perhaps still is the preferred method for many differential geometers. Comparison theorems based on estimates of Jacobi fields and the cut-locus distance were very popular in those days. A key ingredient is the triangle comparison theorem of Toponogov under sectional curvature bounds. This led to an extensive modern research area in ``covariant geometry'' known as metric geometry dealing with length spaces (and other ``limit spaces'') that satisfy curvature bounds from above (mostly for sectional curvature) or from below (sectional and Ricci curvature). Those methods work very well for understanding spaces with bounds for the sectional and Ricci curvature but usually fail to deal with scalar curvature. Schoen and Yau were the first, I believe, who started using the second variation formula for minimal surfaces to deal with problems related to scalar curvature~\cite{SY1}.

\subsection{Jacobi fields and the Ricatti equation}

Here is a simple derivation of the equation for Jacobi fields. If $c(s,t)$ is a variation of geodesics, where $t$ is the arclength parameter of a family of geodesics parametrised by $s$, then $\nabla_X X = 0 $ and $ \nabla_X Y - \nabla_Y X = [X , Y] = 0 $, where $X = \frac{\partial c}{\partial t} , Y = \frac{\partial c}{\partial s}$. Therefore
\[
\nabla_X\nabla_X Y = \nabla_X\nabla_Y X = R(X,Y)X + \nabla_Y\nabla_X X = R(X,Y)X,
\]
which is the second order linear differential equation for a Jacobi field~$Y$ and describes the derivative of the exponential map in terms of curvature.

It is remarkable that from this simple equation and the integrated index form associated with it, one can derive, not just the Bonnet--Meyers theorem, but also a lot of results in Differential Geometry that I learned during the 1970s: from ``pinching" theorems to Bott periodicity (which I learned from Milnor's book ``Morse theory''). Needless to say that there have been various extensions and refinements of this basic equation. For example, to deal with volume estimates in terms of Ricci curvature, Gromov found a crucial generalisation of an original comparison theorem by Bishop.

It turns out that instead of the Jacobi equation, it is much more effective to use the corresponding first order non-linear Riccati equation for the variation of the shape operator (second fundamental form) of a family of hypersurfaces. The derivation is equally simple. If $A = - \nabla N$ denotes the shape operator of a variation of hypersurfaces with unit normal vector $N$ so that $\nabla_N N = 0 $ then
\begin{gather*}
(\nabla_N A)(Y) = - \nabla_N \nabla_Y N + \nabla_{\nabla_N Y} N = R(Y,N)N + \nabla_{[Y , N]} N + \nabla_{\nabla_N Y} N \\
\hphantom{(\nabla_N A)(Y)}{} = R(Y,N)N + A^2(Y).
\end{gather*}

The second variation formula for volumes and areas follow from this by looking at determinants and taking traces. The formula is particularly simple in the case of minimal surfaces in a $3$-dimensional ambient space. The main problem with using minimal surfaces in higher dimensions is the occurrence of possible singularities, but that can be overcome as is shown recently by several researchers (see, for example,~\cite{Lo2, SY5}).

\subsection{Contravariant methods}

This is more about differential forms, spinors, Laplacians, Dirac operators, etc. Exterior differential calculus was a favourite subject of Elie Cartan and S.S.~Chern. Personally, I still like connections and curvatures defined by differential forms on bundles because I learnt Chern--Weil theory of characteristic classes and K-theory at about the same time that I learnt Jacobi fields! Trying to understand the heat kernel proof of the Atiyah--Singer index theorem led me to spinors, which I still find very mysterious. I find it surprising that Hopf (and even Chern) never discussed spinors. Spinors are objects that are more sensitive (or shall I say half as sensitive?) to the action of the orthogonal group than ordinary vectors and exterior forms. In fact they are ``square roots'' of differential forms (the tensor product of spinors is the exterior algebra). This was also one of the reasons why Dirac introduced his operator in quantum physics. The simplest example of spinors is the Hopf bundle: $S^3$ over $S^2$. The bundle is half as curved as the tangent bundle of~$S^2$, both geometrically and topologically. Parallel translation in this bundle around a~great circle in~$S^2$ rotates a~spinor by an angle of $\pi$ instead of $2\pi$ for a~vector and the Euler characteristic of the Hopf bundle is $1$, which is half that of the tangent bundle. More generally, the existence of such a~double cover of the orthonormal frame bundle called a spin structure would be guaranteed by the vanishing of a suitable characteristic class (the second Stiefel--Whitney class). Even if a manifold does not admit a spin structure, twisted spin bundles can still exist, as in the case of K\"ahler manifolds (twisted with the square root of the canonical bundle, for example). The point is that even though the spinor bundle~${\mathbb S}$ and the bundle~$E$ may not exist globally on a manifold~${\mathbb S} \otimes E$ can be well-defined if the ambiguities (or the many-valued-ness) of the two bundles just cancel each other out!

\subsection{The Lichnerowicz formula and the Atiyah--Singer index theorem}

The Bochner technique of proving vanishing theorems for harmonic forms rely on expressing the relevant Laplacian as
a sum of a non-negative operator (the rough Laplacian) and a purely algebraic terms depending only on the curvature.
For the square of the Dirac operator $D$ acting on spinors the corresponding result is the famous formula of Lichnerowicz
which he proved very shortly after the publication of the Atiyah--Singer index theorem:
\[
D^2 = \nabla^\ast\nabla + \frac{R}{4},
\]
where $\nabla$ is the Levi-Civita connection, $\nabla^*$ its adjoint and $R$ is the scalar curvature.

The surprising thing here is the simplicity of the curvature term. Only the scalar curvature appears. As we will see in the proof below, this is partly because the spin representation is very ``democratic'' in the sense that all weights are equal.

\begin{proof} We first define the covariant Hessian
\[ \nabla^2_{u,v} = \nabla_u\nabla_v - \nabla_{\nabla_{u}v},\]
$\nabla^2_{u,v}$ is tensorial in $u$, $v$ and its antisymmetric part is the curvature:
$R(u,v) = \nabla^2_{u,v} - \nabla^2_{v,u}$ (here we use the fact that the Levi-Civita connection is torsion free).
Using a frame satisfying $\nabla_{e_k} e_l =0$ at a given point we have $D = \sum_k e_k \cdot \nabla_{e_k}$, where the dot is Clifford multiplication, and hence
\begin{gather*}
D^2 = \sum_{k,l} e_k \cdot \nabla_{e_k} \bigl( e_l \cdot \nabla_{e_l}\bigr)
 = \sum_{k,l} e_k \cdot e_l \cdot \nabla_{e_k} \nabla_{e_l} \\
\hphantom{D^2}{} = \sum_{k=l} e_k \cdot e_l \nabla_{e_k} \nabla_{e_l} +
2 \sum_{k<l} e_k \cdot e_l \nabla_{e_k} \nabla_{e_l}
= -\sum_{k} \nabla_{e_k} \nabla_{e_k} + \sum_{k<l} e_k \cdot e_l \cdot R(e_k,e_l).
\end{gather*}

The first term is the rough Laplacian $\nabla^* \nabla$ and the second term can be simplified as
\[
 \sum_{k<l} e_k \cdot e_l \cdot R(e_k,e_l)
= - \sum^{m}_{a=1} e_a \cdot \hat{R}(e_a),
\]
where $\{ e_a\}$, $a=1,\dots,m=\frac{n(n-1)}{2}$ is now an orthonormal base for $\wedge ^2 (TM)$ and
$\hat{R}$ is the curvature operator of the Riemannian manifold (note the sign change).

Now choose a base $\{ e_a\}$ for $\wedge^2 (TM)$ that diagonalizes the curvature operator: $\hat{R}(e_a)=\lambda_a e_a $. This acts on a spinor $\psi$ like $\hat{R}(e_a)\psi = \frac{1}{2}\lambda_a e_a \cdot \psi$ (note the $\frac{1}{2}$), and so
\[
- \sum^{m}_{a=1} e_a\cdot \hat{R}(e_a)\psi = - \frac{1}{2}\sum^{m}_{a=1} \lambda_a e_a \cdot e_a \cdot \psi=\frac{R}{4} \psi.
\]
 This proves the Lichnerowicz formula. The formula implies that a compact spin manifold with positive scalar curvature has no non-zero harmonic spinors, and so by the index theorem, compact spin manifolds with non-zero $\hat{A}$-genus do not carry metrics of positive scalar curvature.

Integrating the Lichnerwicz formula on a manifold with boundary, we obtain
\[
 \int_M\left(|\nabla\psi|^2 + \frac{R}{4}|\psi|^2\right) + \int_M | D \psi|^2 =
\int_{\partial M}\langle \nabla_\nu\psi + \nu \cdot D\psi,\psi\rangle,
\]
where $\psi$ is a spinor and $\nu $ is the unit outer normal vector of the boundary.

The formula is proved by computing the divergence of a one form and applying Stokes'
theorem. The specific one form $\alpha $ we use here is defined by
\[
\alpha (v) = \langle\nabla_v\psi + v \cdot D\psi,\psi\rangle
\]
for $v \in TM$ and for a fixed spinor field $\psi$. The boundary operator can be also written as
\begin{gather*}
\nabla_\nu + \nu D =
 \nu \widehat D - \frac{H}{2},
\end{gather*}
where $H$ is the mean curvature of the boundary, and $\nu \widehat D$ is a tangential self adjoint boundary operator, which is useful for imposing Atiyah--Patodi--Singer type non-local boundary conditions. For a harmonic spinor, we then have
\begin{eqnarray*}
\int_M\left(|\nabla \psi|^2 + \frac{R}{4} |\psi|^2\right)
 = -\int_{\partial M} \langle \nu \widehat D\psi, \psi \rangle +
\int_{\partial M} \frac{H}{2} |\psi|^2\ .
\end{eqnarray*}

For twisted Dirac operators acting on spinors with values in a vector bundle
$E$, the Lichnerowicz formula for $D^2$ is computed to be
\[
D^2(\psi \otimes \phi) = \nabla^*\nabla(\psi \otimes \phi) + \frac{R}{4}\psi \otimes \phi + {\mathcal R}(\psi \otimes \phi)
\]
for $\psi \otimes \phi \in \Gamma ({\mathbb S} \otimes E)$, where the last extra term is given by
\[
{\mathcal R}(\psi \otimes \phi) = \frac{1}{2}
\sum^{m}_{j,k=1} e_a \cdot \psi
\otimes R^{\nabla}(e_a)\phi,
\]
$\{e_{a}\}$, $a=1,\dots,m=\frac{n(n-1)}{2}$ is an orthonormal for $\wedge^{2}(T_pM)$ and~$R^{\nabla}$ is the curvature of $E$.

This shows that if the scalar curvature is large (and positive) compared to the curvature $R^{\nabla}$ of the twisting bundle $E$, then there would be no harmonic $E$-valued spinors and the index of the twisted Dirac operator has to vanish, so compact manifolds which allow almost flat twisting bundles which have non vanishing index cannot allow positive scalar curvature.

These twisted Dirac operators play a fundamental role ($K$-theoretically and otherwise) in the Atiyah--Singer index theorem (see~\cite{BGV}) and their index is
\[\int_{M} \hat A(M) \wedge \operatorname{ch}(E).\]

It is interesting to note that in the Seiberg--Witten equations for $4$-manifolds where $E$ is a line bundle, the harmonic spinor~$\psi$ is coupled to the curvature~$R^{\nabla}$, (more precisely to the self-dual part of it) via a natural quadratic form (spinors are square roots of differential forms after all, so if you square them you can get a $2$-form!).
\end{proof}

\section{Scalar curvature}

\subsection{Scalar curvature on the torus}

A torus cannot carry a metric of positive scalar curvature. In fact, any metric of non-negative scalar curvature on a torus is flat. Gromov and Lawson \cite{GL1,GL2,GL3} proved that using the Lichnerowicz formula for twisted Dirac operators and the index theorem. It was a contravariant proof! Gromov~\cite{G1} later gave a conceptual explanation of the main idea by introducing the notion of K-area.

A general principle in Riemannian geometry states that large positive curvature should imply ``small size''. For Ricci curvature, this is made precise by the Bonnet--Myers theorem. This is of course no longer true if the Ricci curvature is replaced by scalar curvature, since a Riemannian product with a sufficiently small $S^2$ has arbitrarily large positive scalar curvature. One might therefor expect that a manifold with large positive scalar curvature is small in the sense that it is ``close'' to a codimension~$2$ subvariety. Gromov~\cite{G1} introduced the notion of $K$-area to quantify this intuition. The $K$-area is, roughly speaking, the inverse of the norm of the smallest curvature obtainable among all topologically essential unitary bundles equipped with connections on a given Riemannian manifold. To measure the norm of the curvature, the metric $g$ is used. However, the definition does not involve the Riemannian curvature of the metric and the $K$-area is a pure $C^0$-invariant of $g$, or more precisely of the metric on 2-forms. It measures the $K$-theoretic 2-dimensional size of the manifold.

\begin{Definition}
\[
K\text{-area} \big(M^{2m},g\big) = \sup_{E,\nabla}\big\| R^\nabla (E)\big\|^{-1},
\]
where $\|~\|$ is the maximum norm and the supremum is taken over all homologically essential unitary bundles
$E$ of all dimensions and over all linear connections~$\nabla$. Homologically essential is equivalent to the fact that the twisted Dirac operator on $E$-valued spinors has a non-zero index, i.e., $\int_{M} \hat A(M) \wedge \operatorname{ch}(E) \neq 0 $.
\end{Definition}

(To extend the definition to odd dimensional manifolds one first defines the K-area for non compact even dimensional manifold using compactly supported bundles~$E$ and characteristic classes. Then one stabilises odd dimensional manifolds by taking products with ${\mathbb R}^{2k+1}$ (see~\cite{G1}).)

The fundamental $K$-area inequality of Gromov can now be stated as follows:

\begin{Theorem} Every complete Riemannian spin manifold $\big(M^n,g\big)$ with scalar curvature \linebreak \mbox{$R(g)
\geq \kappa^2$} everywhere satisfies
\[
K \hbox{\rm -area} (M,g) \leq \frac {c(n)} {\kappa^2}
\]
for some universal constant $c(n)$ depending only on the dimension.
\end{Theorem}

The proof is an immediate consequence of the Lichnerowicz formula, the index
theorem and the definitions. It is an easy observation (by taking large coverings) that the torus has infinite K-area. There has been more precise results about situations when this inequality is sharp, following Llarull's~\cite{Ll1} result for the round sphere. It would also be interesting to find out whether this is related to invariants in symplectic geometry.

A perhaps unrelated, but an important question, is about Einstein metrics. I don't even know whether a torus can carry an Einstein metric with negative scalar curvature.

\begin{Conjecture}
The only Einstein metrics on compact nilmanifolds are the flat metrics on tori.
\end{Conjecture}

\subsection{Positive mass theorems}

Many of the classical concepts of physics, such as mass, energy and momentum are ill-defined in Einstein's theory of general relativity. There is no satisfactory notion of total energy, since the energy of the gravitational field is described in terms of geometry and does not contribute directly to the local stress-energy-momentum tensor. The metric itself is a dynamic variable.
A~measurement is always ``relative'' and hence one has to ``break the symmetry'' in order to have a reasonable notion of mass and energy. The most natural situation occurs when gravitational forces are weak everywhere, except for a confined isolated region. More precisely, in an asymptotically flat space time describing an isolated system like a star or a black hole, where the gravitational field approaches ordinary Newtonian gravity at infinity, one can define the total mass by comparison with Newtonian theory at large distances, just as in the model case of the Schwarzschild metric. Although Hermann Weyl was the first to propose a tentative definition of the energy and mass of an isolated
system in his book: ``Space, time and matter'', it was Arnowitt, Deser and Misner who later gave a more precise definition.
The positive mass theorem states that an asymptotically flat space-time satisfying Einstein's field equations where the stress energy-momentum tensor
satisfies the dominant energy condition has non-negative mass and the mass is zero if and only if it is the flat metric. By restricting to an asymptotically Euclidean space-like slice that is used to define the mass, it can be translated to a problem in Riemannian geometry and the statement becomes:

\begin{Theorem}
An asymptotically Euclidean $3$-manifold $M$ with non-negative scalar
curvature everywhere has positive total mass $m$. Moreover $m=0$ if and only if
$M$ is isometric to flat Euclidean space.
\end{Theorem}

Asymptotically Euclidean means that with respect to appropriate coordinates
\[
 \partial_{\alpha} ( g_{ij} - \delta_{ij}) \in O\big(r^{-1 - |\alpha|}\big)\qquad \text{for} \quad |\alpha| \leq
 2
\]
 and
\[
 m = \frac{1}{16\pi}\lim_{r \rightarrow \infty}{\oint }_{S(r)} (\partial_kg_{ik} - \partial_ig_{kk}) \,{\rm d}\sigma^i.
\]
Here $S(r)$ is the sphere of ``radius'' $r$ in those asymptotic coordinates. Both Schoen--Yau~\mbox{\cite{SY3,SY4}} and Witten~\cite{W1} established this theorem in the course of their proof of the positive mass conjecture. For general dimensions, this was first proved by Bartnik~\cite{Ba}. For a compactly supported perturbation of the flat metric, the above theorem would be a simple
consequence of the result for the torus.

Witten's proof of the positive mass conjecture was contravariant. Motivated by supersymmetry, he uses spinors and a Lichnerowicz type formula with an asymptotic boundary term which can then be identified with the mass at infinity. The sign of the curvature term in the bulk is controlled by Einstein's equations and a positivity condition (the dominant energy condition) for the stress energy momentum tensor. The assumption on the asymptotic geometry gives rise to a trivialisation by parallel spinors (asymptotic infinitesimal supersymmetries) at infinity which are then extended to the whole space-like slice to satisfy an elliptic equation, namely an appropriate Dirac equation. This is the main analytical step. The rest follows (as usual!) by the Lichnerowicz formula with a boundary term. In fact, there is even no need to use the index theorem. Stokes' theorem suffices!

Since Ricci flat metrics are crucial not just in General Relativity but also in many areas of modern theoretical physics (string theory, supergravity etc.), it would be interesting to study the following:

\begin{Problem} Investigate scalar curvature rigidity of asymptotically Ricci flat metrics.
\end{Problem}

Perhaps I should mention that the only time I had a chance to speak with Witten (for a few minutes in 1985), he actually told me to look at that problem! Most of the Ricci-flat metrics we know have special holonomy (I include ${\rm SU}(n)$ among the special holonomy groups), so they admit parallel spinors, and if we follow Witten's proof we should extend them harmonically to the whole manifold and then apply the Lichnerowicz formula with an asymptotic boundary term assuming that the scalar curvature is non-negative. In this context it would also be good to know the answer to the following:

\begin{Question} Are there Ricci flat metrics that do not have special holonomy groups?
\end{Question}

\subsection{Scalar curvature on hyperbolic spaces}

During the late 1980s, I was able to modify Witten's ideas to prove the following result about the scalar curvature rigidity of hyperbolic space~\cite{M1}:

\begin{Theorem}
An asymptotically hyperbolic spin manifold of dimension ${>}2$, whose
scalar curvature satisfies $R \geq -n(n-1)$ everywhere, is isometric to
hyperbolic space.
\end{Theorem}

The proof is an adaptation of Witten's proof but it involves using a hyperbolic connection on an extended bundle, instead of the usual Levi-Civita connection on the tangent bundle. The hyperbolic connection has zero curvature for hyperbolic space (and Euclidean space has positive curvature!) I first used this sort of connection, called a Cartan connection, to prove a ``pinching theorem'' for complex projective space in my dissertation of 1976.

In analogy with the relation between the scalar curvature rigidity of the flat torus and the positive mass theorem, one can ask whether the following is true.

\begin{Question}
Let $g$ be a Riemannian metric on a compact quotient of hyperbolic space $M =
\Gamma \backslash H^n$ of dimension $n \geq 3 $ satisfying $\operatorname{vol}(M,g) =
\operatorname{vol}(M,\bar g)$ where $\bar g$ is the hyperbolic metric of constant
sectional curvature $-1$. If the scalar curvature satisfies: $R(g)\geq R(\bar
g)\equiv - n(n-1)$ everywhere, is $g$ then isometric to $\bar g$?
\end{Question}

If we replace the assumption on scalar curvature by Ricci curvature then this is true by a~deep entropy rigidity
result of Besson--Courtois--Gallot~\cite{BCG}, which is already a vast generalisation of Mostow rigidity. In fact, the question is a special case of a much more general conjecture of Gromov.

\begin{Conjecture}[Gromov]
Let $M^n$ be a complete Riemannian manifold with sectional curvature $K \leq 1$ and let $N^n$ be a compact Riemannian manifold with scalar curvature $R \geq -n(n-1)$. Then every continuous map $f_0\colon M \rightarrow N $ is homotopic to a map $f_1$, such that $\operatorname(f_1(M)) \leq \operatorname(N)$. Moreover, this inequality is strict, unless $N$ has a constant negative curvature and the map $f_0$ is homotopic to a locally isometric map.
\end{Conjecture}

\subsection{Scalar curvature on the hemisphere}

During the early 1990s, I optimistically announced the following conjecture about the scalar curvature rigidity of the hemisphere, after unsuccessfully attempting to prove it with spinorial methods using Cartan connections in the spirit of my proof in the hyperbolic case.

\begin{Conjecture}[false]
Let $M^n$ be a compact spin manifold with simply connected boundary~$\partial M$ and
let $g$ be a Riemannian metric on $M$ with the following properties:
\begin{enumerate}\itemsep=0pt
\item[$(i)$] $\partial M$ is totally geodesic in $M$ and $\partial M$ has constant sectional curvature $K \equiv 1$;
\item[$(ii)$] the scalar curvature of $g$ satisfies $R \geq n(n-1)$ everywhere on M.
\end{enumerate}
Then $(M, g)$ is isometric to the round hemisphere with the standard metric.
\end{Conjecture}

This was proven to be false almost 20 years later, in 2011, by S.~Brendle, F.~Marques and A.~Neves~\cite{BMN}. They found an explicit counter-example. However my conjecture is true if one replaces the scalar curvature bound by a stronger bound on the Ricci curvature as was proved by F.~Hang and X.~Wang in 2009~\cite{HW2}. it is also true if one restricts to the conformal class of the standard metric as was shown earlier in 2006~\cite{HW1}. My false conjecture (fortunately?) has led to a series of interesting papers (see~\cite{B} for a review). These papers (unfortunately?) do not use Dirac operators and spinors. It would be intriguing to understand the correct boundary value problems for harmonic spinors that would lead to a better understanding of non-rigidity or rigidity of lower bounds on the scalar curvature on manifolds with prescribed geometry near the boundary. Are the Atiyah--Patodi--Singer boundary conditions the right ones to use?

\subsection{Gromov's recent work on scalar curvature}

In a series of recent papers \cite{G2,G3,G4}, Gromov has been looking at scalar curvature from a~much broader and more geometric point of view, beyond Dirac operators and spinors. He established a~number of new results, including new proofs based on recent regularity results for minimal surface and soap bubbles in higher dimensions established by various researchers, for example~\mbox{\cite{Lo2, SY5}}. An apparent advantage of the covariant method is the fact that it requires less smoothness, so it can deal with ``rougher'' metrics and even singularities. In fact, one of Gromov's goals is to have a more discrete and combinatorial understanding of scalar curvature for polyhedral objects. Gromov does lament the fact that one still does not have a good understanding of how the two methods that I have described in this essay relate to each other. I~should mention that the late S.S.~Chern also made that remark to me in 1988 after I gave a talk about my extension of Witten's proof to the hyperbolic case. He told me that it was my ``homework" to relate Schoen--Yau's covariant proof to Witten's contravariant proof. I am still working on it!

\section{Some speculative contravariant remarks}

As should be obvious by now, I am somewhat biased towards the ``softer'' (simpler?) contravariant methods, so in conclusion, let me make some vague and speculative remarks about using spinors and the index theorem in differential geometry.
\begin{enumerate}\itemsep=0pt
\item In proving sharp K-area estimates in terms of a positive lower bound for the scalar curvature, the main trick is to find the optimal twisting bundle $E$ with non-zero index (i.e., homologically essential) and the optimal connection $\nabla$ on $E$ that is just right for the term involving the curvature of the bundle~$R^{\nabla}$ to cancel the scalar curvature term in the Lichnerowicz formula. Is there a general method to determine those bundles and connections that optimise Gromov's $K$-area inequality? Can we detect them in the classifying space of these bundles? More ambitiously can we couple the harmonic spinor to the curvature of the bundle as is done in the Seiberg--Witten equations to really balance the two curvature terms in the Lichnerowicz formula? These ideas are also useful for understanding $K$-length~\cite{G1} which gives information about spectral gaps.

\item Almost all ``vanishing theorems'', such as that of Lichnerowicz or Bochner, prove the vanishing of all harmonic objects. This is, at least superficially, a lot stronger than just saying that some index of an elliptic operator is zero. The index is the asymmetry or disparity between harmonic objects of different parities (an anomaly, in physical lingo). A~vanishing index only says that the two types of harmonic objects (zero modes) are of the same dimension. It doesn't say, like a vanishing theorem, that there are no zero modes at all. I find it intriguing that in many cases (BPS-states?) one type of zero modes is already automatically excluded, so the index is actually the dimension of harmonic objects of a~certain type.

\item Is there a direct way to write down explicitly the closed characteristic differential form representing the $\hat {A}$-genus of a compact closed spin $4k$-dimensional manifold with positive scalar curvature as an exact form ${\rm d} \eta$? One thing about differential forms representing Pontrjagin classes is that they depend only on the conformally invariant Weyl curvature tensor which a priori has nothing to do with the scalar curvature! The index density is the (super-)trace of the heat kernel $\exp\big({-}tD^2\big)$ as $t \to 0^+$ whereas the harmonic projectors describe the behaviour of the heat kernel as $t \to \infty$, so maybe there is a way to relate the index density and the Lichnerowicz formula by varying~$t$ in the heat kernel?

\item The Euler characteristic can be localized around zeros of vector fields. Is there any way of localizing the $\hat A$-genus around codimension~2 submanifolds, relating it perhaps to the scalar curvature and the mean curvature?

\item The index theorem is best proved by evaluating the supertrace of the heat kernel of the elliptic operator on the diagonal. Is there a version of the index theorem where the heat kernel is evaluated differently, for example outside of the diagonal, by letting one point go to infinity at a specific rate with respect to~$t$. One might not be interested in the classical topological index, but rather in the index (probability) density at the infinite (or finite) boundary. This is probably related to the Callias-type index theorem and perhaps to the AdS/CFT correspondence (holographic principle).

\item Another speculation is that the local asymptotic expansions (as $t \to 0^+$) of the heat kernel are always done with respect to the flat Euclidean metric as the background geometry. Curvature is always measured with respect to Euclidean geometry, which by definition is flat. My use of Cartan connections show that curvature can be measured with respect to other background metrics, especially that of symmetric spaces. It would be interesting to find an ``index theorem" where the asymptotics of the heat kernel are calculated with respect to other background geometries, both locally and at infinity. The result might be more geometrical (in the spirit of AdS/CFT) than a purely topological statement about the index of an elliptic operator. One should also extend the Bochner--Weitzenb\"ock formulas as I have done for the Lichnerowicz formula in the hyperbolic case.

\item To study deformations and obstructions to Einstein metrics, I suggest that one should look at ``secondary characteristics classes'' and other cycles, (coming from characteristic classes?) on the Lie algebra of vector fields, i.e., the Lie algebra of the diffeomorphism group, which acts on the space of metrics. This might be useful to study Einstein metrics with negative scalar curvature. One should also study the Rarita--Schwinger operator for $\frac{3}{2}$-spinors (which describe gravitinos in physics), instead of just the Dirac operator for $\frac{1}{2}$ spinors (which are just infinitesimal supersymmetries, or square roots of vector fields after all!).
\end{enumerate}

\pdfbookmark[1]{References}{ref}
\LastPageEnding


\begin{thebibliography}{99}
\footnotesize\itemsep=0pt

\bibitem{Ba}
Bartnik R., The mass of an asymptotically flat manifold, \href{https://doi.org/10.1002/cpa.3160390505}{\textit{Comm. Pure
 Appl. Math.}} \textbf{39} (1986), 661--693.

\bibitem{BGV}
Berline N., Getzler E., Vergne M., Heat kernels and {D}irac operators,
 \textit{Grundlehren der Mathematischen Wissenschaften}, Vol.~298, \href{https://doi.org/10.1007/978-3-642-58088-8}{Springer-Verlag}, Berlin, 1992.

\bibitem{BCG}
Besson G., Courtois G., Gallot S., Entropies et rigidit\'es des espaces
 localement sym\'etriques de courbure strictement n\'egative, \href{https://doi.org/10.1007/BF01897050}{\textit{Geom.
 Funct. Anal.}} \textbf{5} (1995), 731--799.

\bibitem{B}
Brendle S., Rigidity phenomena involving scalar curvature, \href{https://doi.org/10.4310/SDG.2012.v17.n1.a4}{Int. Press}, Boston,
 MA, 2012, 179--202, \href{https://arxiv.org/abs/1008.3097}{arXiv:1008.3097}.

\bibitem{BMN}
Brendle S., Marques F.C., Neves A., Deformations of the hemisphere that
 increase scalar curvature, \href{https://doi.org/10.1007/s00222-010-0305-4}{\textit{Invent. Math.}} \textbf{185} (2011),
 175--197, \href{https://arxiv.org/abs/1004.3088}{arXiv:1004.3088}.

\bibitem{G1}
Gromov M., Positive curvature, macroscopic dimension, spectral gaps and higher
 signatures, in Functional Analysis on the Eve of the 21st Century,
 {V}ol.~{II} ({N}ew {B}runswick, {NJ}, 1993), \href{https://doi.org/10.1007/s10107-010-0354-x}{\textit{Progr. Math.}}, Vol.~132,
 Birkh\"auser Boston, Boston, MA, 1996, 1--213.

\bibitem{G2}
Gromov M., Dirac and {P}lateau billiards in domains with corners, \href{https://doi.org/10.2478/s11533-013-0399-1}{\textit{Cent.
 Eur.~J. Math.}} \textbf{12} (2014), 1109--1156, \href{https://arxiv.org/abs/1811.04318}{arXiv:1811.04318}.

\bibitem{G3}
Gromov M., Metric inequalities with scalar curvature, \href{https://doi.org/10.1007/s00039-018-0453-z}{\textit{Geom. Funct.
 Anal.}} \textbf{28} (2018), 645--726, \href{https://arxiv.org/abs/1710.04655}{arXiv:1710.04655}.

\bibitem{G4}
Gromov M., Scalar curvature of manifolds with boundaries: natural questions and
 artificial constructions, \href{https://arxiv.org/abs/1811.04311}{arXiv:1811.04311}.

\bibitem{GL1}
Gromov M., Lawson Jr. H.B., Spin and scalar curvature in the presence of a
 fundamental group.~{I}, \href{https://doi.org/10.2307/1971198}{\textit{Ann. of Math.}} \textbf{111} (1980), 209--230.

\bibitem{GL2}
Gromov M., Lawson Jr. H.B., The classification of simply connected manifolds of
 positive scalar curvature, \href{https://doi.org/10.2307/1971103}{\textit{Ann. of Math.}} \textbf{111} (1980),
 423--434.

\bibitem{GL3}
Gromov M., Lawson Jr. H.B., Positive scalar curvature and the {D}irac operator
 on complete {R}iemannian manifolds, \href{https://doi.org/10.1007/BF02953774}{\textit{Inst. Hautes \'Etudes Sci. Publ.
 Math.}} \textbf{58} (1983), 83--196.

\bibitem{HW1}
Hang F., Wang X., Rigidity and non-rigidity results on the sphere,
 \href{https://dx.doi.org/10.4310/CAG.2006.v14.n1.a4}{\textit{Comm. Anal. Geom.}} \textbf{14} (2006), 91--106.

\bibitem{HW2}
Hang F., Wang X., Rigidity theorems for compact manifolds with boundary and
 positive {R}icci curvature, \href{https://doi.org/10.1007/s12220-009-9074-y}{\textit{J.~Geom. Anal.}} \textbf{19} (2009),
 628--642, \href{https://arxiv.org/abs/0911.0380}{arXiv:0911.0380}.

\bibitem{Ll1}
Llarull M., Sharp estimates and the {D}irac operator, \href{https://doi.org/10.1007/s002080050136}{\textit{Math. Ann.}}
 \textbf{310} (1998), 55--71.

\bibitem{Lo2}
Lohkamp J., Minimal smoothings of area minimizing cones, \href{https://arxiv.org/abs/1810.03157}{arXiv:1810.03157}.

\bibitem{M1}
Min-Oo M., Scalar curvature rigidity of asymptotically hyperbolic spin
 manifolds, \href{https://doi.org/10.1007/BF01452046}{\textit{Math. Ann.}} \textbf{285} (1989), 527--539.

\bibitem{SY1}
Schoen R., Yau S.-T., Existence of incompressible minimal surfaces and the
 topology of three-dimensional manifolds with nonnegative scalar curvature,
 \href{https://doi.org/10.2307/1971247}{\textit{Ann. of Math.}} \textbf{110} (1979), 127--142.

\bibitem{SY3}
Schoen R., Yau S.-T., On the proof of the positive mass conjecture in general
 relativity, \href{https://doi.org/10.1007/BF01940959}{\textit{Comm. Math. Phys.}} \textbf{65} (1979), 45--76.

\bibitem{SY4}
Schoen R., Yau S.-T., Proof of the positive mass theorem.~{II}, \href{https://doi.org/10.1007/BF01942062}{\textit{Comm.
 Math. Phys.}} \textbf{79} (1981), 231--260.

\bibitem{SY5}
Schoen R., Yau S.-T., Positive scalar curvature and minimal hypersurface
 singularities, \href{https://arxiv.org/abs/1704.05490}{arXiv:1704.05490}.

\bibitem{W1}
Witten E., A new proof of the positive energy theorem, \href{https://doi.org/10.1007/BF01208277}{\textit{Comm. Math.
 Phys.}} \textbf{80} (1981), 381--402.

\end{thebibliography}
\end{document}